\documentclass[a4paper]{amsart}
\usepackage{enumerate}
\usepackage{amssymb}
\usepackage{amsmath}
\usepackage{amsthm}
\usepackage{float}
\usepackage{verbatim}
\usepackage{url}

\def\F{{\mathbb F}}
\newcommand{\sli}{{\mathfrak{sl}}}
\def\ad{{\rm ad}}
\def\End{{\rm End}}

\def\Extr{{\mathcal E}}
\def\Sand{{\mathcal S}}
\def\la{\langle}
\def\ra{\rangle}
\newcommand{\Ker}{\mathrm{Ker}\,}

\newcommand{\pd}{\partial}
\renewcommand{\phi}{\varphi}

\newcommand{\W}{{W_{1,1}(5)}}
\newcommand{\wW}{\widetilde{\W}}
\newcommand{\Z}{\mathbb{Z}}

\newtheorem{theorem}{Theorem}[section]
\newtheorem{lemma}[theorem]{Lemma}
\newtheorem{proposition}[theorem]{Proposition}
\newtheorem{example}[theorem]{Example}
\newtheorem{notation}[theorem]{Notation}

\renewenvironment{proof}[1][Proof]{\begin{trivlist}\item[\hskip \labelsep {\bfseries #1}]}{\qed\end{trivlist}}

\newcounter{storedequation} 
\let\storedtheequation=\theequation
\newenvironment{lettereqns}[2][0]{%
  \setcounter{storedequation}{\value{equation}}%
  \setcounter{equation}{#1}%
  \renewcommand{\theequation}{#2\arabic{equation}}}{%
  \setcounter{equation}{\value{storedequation}}%
  \renewcommand{\theequation}{\storedtheequation}}

\title{Simple Lie Algebras having Extremal Elements}

\author[Arjeh~M.~Cohen]{Arjeh~M.~Cohen}
\address[Arjeh~M.~Cohen]{
Department of Mathematics and Computer Science
Technische Universiteit Eindhoven
P.O. Box 513, 5600 MB Eindhoven, Netherlands}
\email{amc@win.tue.nl}

\author[G.~Ivanyos]{G\'abor~Ivanyos }
\thanks{The second author would like to thank the DIAMANT mathematics cluster and an NWO visitor's grant for supporting a visit to Eindhoven, during which part of the work was carried out}
\address[G\'abor~Ivanyos]{
Informatics Research Laboratory, 
Computer and Automation Institute, Hungarian Academy of Sciences,
L\'agym\'anyosi u. 11,
H-1111,
Budapest,
Hungary}
\email{gabor.ivanyos@sztaki.hu}

\author[D.~Roozemond]{Dan~Roozemond}
\address[Dan~Roozemond]{
Department of Mathematics and Computer Science
Technische Universiteit Eindhoven
P.O. Box 513, 5600 MB Eindhoven, Netherlands}
\email{d.a.roozemond@tue.nl}

\begin{document}

\begin{abstract}
Let $L$ be a simple finite-dimensional Lie algebra of characteristic distinct
from $2$ and from $3$.  Suppose that $L$ contains an extremal element that is
not a sandwich, that is, an element $x$ such that $[x,[x,L]]$ is equal to the
linear span of $x$ in $L$.  In this paper we prove that, with a single
exception, $L$ is generated by extremal elements. The result
is known, at least for most characteristics, but the proofs in the literature
are involved. The current proof closes a gap in a geometric proof that every
simple Lie algebra containing no sandwiches (that is, ad-nilpotent elements of
order $2$) is in fact of classical type.
\end{abstract}

\maketitle

\section{Introduction}
Let $L$ be a Lie algebra over a field $\F$ of characteristic distinct from
$2$.  An element $x \in L$ is said to be \emph{extremal} if $[x,[x,L]]
\subseteq \F x$.  If $[x,[x,L]] = 0$ we say $x$ is a \emph{sandwich}.  By
Premet, \cite{Prem,Prem23}, every finite-dimensional simple Lie algebra over
an algebraically closed field of characteristic distinct from $2$ and $3$ is
known to have an extremal element; see \cite{Tange} for a self-contained proof
in case $p>5$.  If a simple Lie algebra is finite-dimensional and generated by
extremal elements, then it is of classical type.  This fact follows from the
classification of finite-dimensional simple Lie algebras as described in
\cite{PSi,PS,Strade}, but can also be derived from geometric arguments using
the theory of buildings, cf. \cite{CIG,CIG2}, up to small rank cases and the
verification that the building determines a unique Lie algebra generated by
extremal elements up to isomorphism -- a subject of ongoing work.  None of
these extremal elements are sandwiches; see \cite[Remark 9.9]{CSUW}.

In order to use these two results for a revision of the classification of
simple Lie algebras of classical type, the gap between the two has to be
filled. In other words, an elementary proof would be needed of the fact that
a simple Lie algebra over an algebraically closed field of characteristic
distinct from $2$ and $3$ having an extremal element that is not
a sandwich is generated by extremal elements.

Using powerful methods, Benkart \cite[Theorem 3.2]{Benkart} shows that if a
simple Lie algebra over an algebraically closed field of characteristic
$p\ge7$ or $p=0$ contains a nilpotent element of order at most $p-1$ and no
sandwiches, then it is of classical type.
Together with the abovementioned results of Premet, this gives that any
simple Lie algebra over an algebraically closed field of characteristic
$0$ or greater than $5$ without sandwiches is of classical type.
Since Benkart's methods are rather
involved, this paper is devoted to a self-contained proof of the
observed gap and an extension to the case of characteristic $5$.  The
field need not be algebraically closed.  Our extension allows for one more
example of a simple Lie algebra having a non-sandwich extremal element, namely
the 5-dimensional Witt algebra $\W$ over a field of characteristic $5$ (see
Example \ref{ex-Witt} for an explicit description of this Lie algebra).  It is
a counterexample in that it only contains one such element up to scalar
multiples.

\begin{theorem}\label{th:main}
Let $\F$ be a field of characteristic distinct from $2$ and $3$, and 
let $L$ be a simple Lie algebra over $\F$.
Suppose that $L$ contains an extremal element that is not a sandwich.
Then either $\F$ has characteristic $5$ and $L$ is isomorphic to $\W$ 
or $L$ is generated by extremal elements.
\end{theorem}

The counterexample was known to Alexander Premet. We are grateful for fruitful
discussions with him about our work.  We would also like to thank Helmut
Strade for the insight he provided us into the classification of modular Lie
algebras.  

We briefly outline the paper.  In Section \ref{sec:dichotomy} we find that a
Lie algebra containing an extremal element that is not a sandwich either has
more extremal elements or is defined over a field of characteristic $5$ and
has a particular Lie subalgebra.  Elementary proofs for most of the statements
in Section \ref{sec:dichotomy} were known before 1977; see \cite{Benkart} and
references therein. We include proofs here for the sake of completeness.  We
gratefully acknowledge David Wales' contribution in the guise of Proposition
\ref{thm_sl2wales}. 
In Section
\ref{sec:p=5}, we pin down the exceptional case in characteristic $5$, and in
Section \ref{sec:fin} we show that, in the absence of the exceptional case,
there are many more extremal elements and finish the proof of the main
theorem.  The proofs of Lemma \ref{lem_ysq0} and Proposition
\ref{prop_noextWitt} were found by experiments with the GAP computer system
package GBNP; see \cite{GBNP}.

To finish the introduction, we fix some notation of use throughout the paper.

\begin{notation} \rm
Throughout this paper,
$\F$ will be a field whose
characteristic is denoted by $p$, and $L$ will be a Lie algebra over $\F$.

An element $x \in L$ is said to be \emph{extremal on $M$} if $[x,[x,M]]
\subseteq \F x$.  If $x$ is extremal on $L$ and no confusion is imminent, we
call $x$ \emph{extremal}.  We write $\Extr_L(M)$ for the set of elements
extremal on $M$ and abbreviate $\Extr_L(L)$ to $\Extr_L$. Furthermore, we
write $\Extr_L(y)$ for $\Extr_L(\{y\})$.

Similarly, if $[x,[x,M]] = 0$ we say $x$ is a \emph{sandwich on $M$},
and if $[x,[x,L]] = 0$ we simply call $x$ a \emph{sandwich}. We write
$\Sand_L(M)$ for the set of sandwiches on $M$. Again, we write
$\Sand_L(y)$ for $\Sand_L(\{y\})$.

By linearity of the expression $[x,[x,m]]$ in $m$, we have $\Extr_L(M)
= \Extr_L(\la M\ra)$, where $\la M\ra$ denotes the linear subspace of
$L$ spanned by $M$.  Hence, when writing $\Extr_L(M)$, we may assume that
$M$ is a linear subspace of $L$, and similarly for $\Sand_L(M)$.
For $x \in \Extr_L(M)$ and $m \in M$, we define $f_x(m) \in \F$ to
be such that $[x,[x,m]] = f_x(m)x$. 
\end{notation}

\section{Jacobson-Morozov type results}\label{sec:dichotomy}
For extremal elements we present a slightly better version of the well-known
theorem by Jacobson and Morozov.  The orginal result, ascribed to Morozov in
\cite[page 98]{Jacobson}, is adapted by David Wales to extremal elements and
works for characteristic at least 5. 

\begin{proposition}\label{thm_sl2wales}
Suppose that $p$ is distinct from $2$ and $3$
and that $L$
contains an extremal element $x$.  
If $w$ is an element for which $f_x(w)=-2$, then, with $h=[x,w]$,
there is $y\in L$ for which 
\begin{eqnarray}
\label{eqn_sl2}
[x,y]=h, \qquad [h,x]=2x,\ \  \mbox{\rm and }\ \ [h,y]=-2y. 
\end{eqnarray}
\end{proposition}

\medskip
The three elements $x$, $y$, $h$ are the usual generators of the Lie algebra
$\sli_2(\F)$ of $2\times2$ matrices of trace $0$ over $\F$.  A triple
satisfying the relations (\ref{eqn_sl2}) is called an $\sli_2$-triple.

\begin{proof}
Let $X=\ad_x$.  Let $h=[x,w]$ and $H=\ad_h$.  The hypothesis
$f_x(w)=-2$  means $[x,[x,w]]=-2x$.  In particular 
$[h,x]=-[x,[x,w]]=2x$ as required.

We know $X$ is nilpotent as $[x,[x,[x,y]]]=[x,f_x(y)x]=0$ and 
so $X^3=0$.  

Let $C_L(x) = \{ u \in L \mid [u,x]=0\} = \Ker X$.
The following computation shows that $[w,h]-2w \in C_L(x)$:
\begin{eqnarray*}
[[w,h]-2w,x] &=& [[w,h],x]-2[w,x]  
= [[w,x],h]+[w,[h,x]]-2(-h) 
\\
&=& [-h,h]+[w,2x]+2h 
= 2(-h)+2h \\
&=&0.
\end{eqnarray*}

Consequently, $[w,h]=2w+x_1$ where $x_1\in C_L(x)$.  
We claim $C_L(x)$ is $H$-invariant.  To see this notice 
$[X,H]=-2X $ and so $XH-HX=-2X$.  Suppose $u\in C_L(x)$.  Then 
$XHu=HXu-2Xu=0$ and so $X(Hu)=0$, proving $Hu \in C_L(x)$ as claimed.

Next we consider the action of $H$ on $C_L(x)$.  Let $u\in C_L(x)$.  
Now, with $W = \ad_w$,
\[
Hu=[X,W]u=XWu-WXu=XWu \in XL,
\]
and 
\begin{eqnarray*}
[X^2,W] &=& X^2W-WX^2 
= X[X,W]+[X,W]X \\
 &=& XH+HX 
= HX-2X+HX = 2(H-1)X,
\end{eqnarray*}
so
\[
2(H-1)Hu=2(H-1)XWu=(X^2W-WX^2)Wu.
\]
But $XW=[X,W]+WX=H+WX$ and so $X^2Wu=XHu+XWXu=0$.  
In particular,
$2(H-1)Hu=X^2W^2u$.
As $x$ is extremal, for any $z \in L$ we have $X^2z=f_x(z)x$, and so 
$
(H-2)X^2z=(H-2)f_x(z)x=0
$
as $Hx=2x$.    
Now 
\[
2(H-2)(H-1)Hu=(H-2)X^2W^2u=0.
\]
This means the eigenvalues of $H$ acting on $C_L(x)$ are $0,1,2$ and 
as the characteristic is at least $5$ we see $-2$ is not 
an eigenvalue. In particular $H+2$ is nonsingular on $C_L(x)$.  
Pick $w_1\in C_L(x)$ for which $(H+2)w_1=x_1$ and so $[h,w_1]=x_1-2w_1$.  
Set $y=w+w_1$.  
Now $[x,y]=[x,w+w_1]=h+0=h$ and so $[x,y]=h$.  Also 
$[h,y]=[h,w+w_1]=(-2w-x_1)+(x_1-2w_1)=-2(w+w_1)=-2y$.  
This completes the proof of the proposition.  
\end{proof}

In the remainder of this section we suppose that $x,y,h \in L$ are an
$\sli_2$-triple.
In fact,
the triple is determined by the pair $x,y$ and the relations
\begin{equation}\label{eqn_sl2_rels}
[[x,y],x]=2x \mbox{ and } [[x,y],y] = -2y,
\end{equation}
as $h = [x,y]$. Such a pair will be called an $\sli_2$-pair.
Note that, if $x$ is extremal, this implies $f_x(y) = -2$.

\begin{lemma}\label{lem_ysq0}
Suppose that $\F$ is of characteristic $p\ne 2,3$ and that $x$ and $y$
are an $\sli_2$-pair in $L$. Set $S = \langle x,y,[x,y] \rangle$.
If $x \in \Extr_L$, then $y$ acts quadratically
on $L/S$, i.e., $\ad_y^2(L/S) = 0$.
\end{lemma}

\begin{proof} \begin{lettereqns}{R}
$S$ is a Lie subalgebra of $L$
isomorphic to $\sli_2(\F)$.  
Consider $L$ as a module on which $S$ acts.
Obviously $S$ is an invariant subspace, so $L/S$ is an $S$-module.  Write $X$,
$Y$ for the action of $\ad_x$, $\ad_y$, respectively, on $L/S$.  As
$\ad_x^2(L) \subseteq \F x \subseteq S$, we have $X^2 = 0$.  We list the
relations (\ref{eqn_sl2_rels}) in terms of $X$ and $Y$, and the quadraticity
of $X$ that we just found.
\begin{eqnarray}
X^2 Y - 2XYX + YX^2 + 2X &=& 0   \\
-XY^2 + 2YXY - Y^2X - 2Y &=& 0   \\
X^2 &=& 0 
\end{eqnarray} 
The relations (R1) and (R3) immediately imply
\begin{eqnarray}
XYX - X = 0.
\end{eqnarray}
Multiplying (R2) from the left by $X$ gives 
$$-X^2Y^2 + 2XYXY - XY^2X -2XY = 0,$$ 
which, after application of (R3) and (R4), gives
\begin{eqnarray}
XY^2X = 0.
\end{eqnarray}
Denote by $R_2$ the left hand side of (R2). Then, by (R3),
\begin{eqnarray*}
0 & = & 
Y R_2 Y X - YXY R_2 + 2 Y^2 X R_2 - R_2 YXY + XYR_2 Y -3YR_2 \cr
& & -2 YXR_2Y + 3 R_2 Y -2YXR_2Y -6 R_2Y + 2XR_2Y^2 \cr
& = & 12Y^2-3XY^3+7YXY^2-5Y^2XY+Y^3X+3XYXY^3\cr
& & -7YXYXY^2+5Y^2XYXY-Y^3XYX.
\end{eqnarray*}
Replacing $XYX$ by $X$ and $X^2$ by $0$, using (R4) and (R3), we find
\begin{eqnarray*}
0 & = & 12Y^2-3XY^3+7YXY^2-5Y^2XY+Y^3X+3XY^3\cr
& & -7YXY^2+5Y^2XY-Y^3X \cr
& = & 12Y^2.
\end{eqnarray*}
As $p\neq2,3$, we conclude that $Y^2=0$.
\end{lettereqns} 
\end{proof} 

For $a\in\End(L)$ and $\lambda\in\F$, we denote by
$L_\lambda(a)$ the $\lambda$-eigenspace of $a$ in $L$.

\begin{theorem}\label{thm_sl2quad}
Suppose that $\F$ is a field of characteristic $p\ne2,3$,
that $L$ is a Lie algebra over $\F$,
and that $x$ is an extremal element of $L$ that is not a sandwich.
Then there are $y$, $h\in L$ such that
$x$, $y$, $h$ is an $\sli_2$-triple.
Moreover, for each such a triple,
$\ad_h$ is diagonizable with eigenvalues $0$, $\pm1$, $\pm2$ and satisfies
$L_{-2}(-\ad_h) = \F x$ and $L_{2}(-\ad_h) = \F y$.
\end{theorem}

\begin{proof} 
As $x$ is not a sandwich and the characteristic of $\F$ is not $2$, there is
$w\in L$ with $f_x(w) = -2$.  By Proposition \ref{thm_sl2wales} with $h =
[x,w]$, there is $y\in L$ such that $x$, $y$, $h$ are an $\sli_2$-triple.
They generate a Lie subalgebra $S$ of $L$ isomorphic to $\sli_2(\F)$.  Viewing
$L$ as an $S$-module as in the proof of Lemma \ref{lem_ysq0}, 
we see that $S$ itself is an invariant submodule.
Denote by $X$, $Y$, and $H$ the actions of $\ad_x$, $\ad_y$, and $\ad_h$,
respectively, on the quotient module $L/S$.  As $x\in\Extr_L$, we have $X^2 =
0$. By Lemma \ref{lem_ysq0}, also $Y^2 = 0$. It readily follows that the
subalgebra of $\End(L/S)$ generated by $X$, $Y$, and $H$ is linearly spanned
by $1$, $X$, $Y$, $H$, $XY$, $XH$, and $YH$, and that the relation $H^3 =H$ is
satisfied. In particular, $H$ is diagonizable on $L/S$ with eigenvalues $0$,
$1$, and $-1$ only.  Consequently, there are subspaces $U$ and $V$ of $L$ such
that $L = S + U + V$ is a direct sum of subspaces such that $(S+U)/S=\Ker
(H^2-1)$ and $(S+V)/S = \Ker H$. Notice that $\ad_h$ has eigenvalues $-2$,
$0$, $2$ on $S$, each with multiplicity $1$. A small computation shows that
actually $\Ker \ad_h^2 = \Ker\ad_h$, so that $-\ad_h$ is diagonizable with
eigenspaces $L_i(-\ad_h)$ $(i=-2,-1,0,1,2)$ satisfying $L_{-2}(-\ad_h) = \F x$
and $L_{2}(-\ad_h) = \F y$.
\end{proof}

To end this section, we exploit the $\ad_h$-grading with five components.
The following result is a slight variation of \cite[Proposition 22]{CIG}.

\begin{proposition}\label{prop_Zgrading}
Suppose $x \in \Extr_L$ and $y \in L$ are an $\sli_2$-pair
in a Lie algebra $L$ of characteristic $p>3$. 
Let $L_i = L_i(-\ad_h)$ $(i=-2,-1,0,1,2)$
be the
components of the $\Z_p$-grading by 
$h = [x,y]$.
Then
either $p=5$ and $[y,[y,v]] = x$ for some $v\in L_{-1}$, or 
$y$ is extremal in $L$, the
components $L_i$ $(i=-2,-1,0,1,2)$ 
actually give a $\Z$-grading of $L$,
with $L_{-2}=\F x$, $L_{2}=\F y$,
$[x,L_{-1}]=L_1$, and $[y,L_{1}]=L_{-1}$.
\end{proposition}

\begin{proof}
Set $S = \F x + \F y + \F h$. By assumption, $S \cong \sli_2(\F)$. 
The identifications of $L_{-2}$ and $L_2$ with $\F x$ and $\F y$, respectively,
were established in Theorem \ref{thm_sl2quad}. Suppose that $y$
is not an extremal element. As $\ad^2_y L_i\subseteq \F y$ for
$i\neq \pm 1$ and $\ad_y L_{-1}\subseteq L_1$, this can only happen
if $[y,L_1]\neq 0$. Then, by the grading properties,
$[y,L_1]\subseteq L_3(-\ad_h)$ and so $3$ is equal to a member $i$ of 
$\{-2,-1,0,1,2\}$ modulo $p$. As $p\ge5$, this implies $p=5$ and
$i= -2$. Thus $[y,L_1]=\F x$. It follows that, for every 
$u\in L_1$, $\ad_x\ad_y u=0$, 
whence $\ad_y \ad_x u=(\ad_x\ad_y-\ad_h) u=-u$. 
Therefore $[y,[y,L_{-1}]\supseteq [y,[y,[x,L_1]]]=[y,L_1]=\F x$,
and, by homogeneity, $[y,[y,L_{-1}]] \subseteq L_{-2} = \F x$,
so the first case holds. 
To complete the proof, assume that both $x$ and $y$ are extremal.
The argument above also shows that if $[y,L_1]\neq 0$
then $p=5$ and $[y,L_1]=\F x$. It follows then
that $[y,[y,L_1]]=\F h\not\subseteq \F y$, a contradiction
to extremality of $y$. Thus $[y,L_1]=0$ and, similarly,
$[x,L_{-1}]=0$. It follows that for every pair $i,j$
from the interval $[-2,2]$, we have
$[L_i,L_j]=0$ whenever the ordinary sum $i+j$ falls
outside the interval $[-2,2]$. Thus the grading is
indeed a $\Z$-grading. To see the very last 
two equalities of the proposition just notice that
for every $u\in L_1$ we have
$\ad_y\ad_x u=-u$ for every $u\in L_1$ and,
similarly, $\ad_x\ad_y v=-v$ for every $v\in L_{-1}$ as observed above.
\end{proof}

\section{The characteristic 5 case}\label{sec:p=5}
Suppose that $p=5$, and that $x$ is an extremal element of $L$ that is not a
sandwich.  By Proposition \ref{thm_sl2wales} there are $y,h \in L$ such that
$x$ and $y$ are an $\sli_2$-triple.  By Theorem \ref{thm_sl2quad}, $\ad_h$ is
diagonizable and there exists a grading of $L$ by $-\ad_h$ eigenspaces $L_i$
$(i=-2,-1,0,1,2)$.  In this section we consider the case where $y$ is not an
extremal element.  By Proposition \ref{prop_Zgrading} there exists an element
$v\in L_{-1}$ such that $[y,[y,v]] = x$.

\begin{example}\label{ex-Witt}
\rm
Before we proceed, we show that this case actually occurs.
The $5$-dimensional Witt algebra $\W$ can be defined as follows.
Let $\F$ be a field of characteristic $p=5$ and take
the vector space over $\F$
with basis $z^i \partial_z$, for $i=0,\ldots,4$. 
The Lie bracket is defined on two of these elements by
\begin{eqnarray*}
[ z^i \pd_z, z^j \pd_z ] & := & z^i \pd_z( z^j \pd_z ) - z^j \pd_z( z^i \pd_z ) 
 = j z^i z^{j-1} \pd_z - i z^j z^{i-1} \pd_z \\
& = & (j-i) z^{i+j-1} \pd_z,
\end{eqnarray*}
with the convention that 
\begin{equation}\label{eqn_witt_conv1}
z^i := 0 \mbox{ whenever } i \notin \{ 0, \ldots, 4\}.
\end{equation}
The Lie bracket extends bilinearly to a multiplication
on $\W$. It is antisymmetric and satisfies
the Jacobi identity, so that $\W$ is indeed a Lie algebra of dimension
$5$ over $\F$.

Now we construct an extension $\widetilde{\W}$ of $\W$: 
Add one basis element, namely $z^6
\pd_z$, and adapt (\ref{eqn_witt_conv1}):
\begin{equation}\label{eqn_witt_conv2}
z^i := 0 \mbox{ whenever } i \notin \{0, 1,2,3,4,6\}.
\end{equation}

The only entry of the multiplication table that differs between $\wW$
and $\W$ is $[z^3 \pd_z, z^4 \pd_z]$: This is $0$ in $\W$ and $z^6
\pd_z$ in $\wW$.  Furthermore, $z^6 \pd_z$ commutes with all other
elements.  So $\wW$ is indeed an extension of $\W$ by a
$1$-dimensional center.  This extension was constructed in
\cite{Block}. The analog over the complex numbers of $\wW$ is also known as the
Virasoro algebra.

Now let $W$ be either $\W$ or $\wW$.
Then $x = -z^2\pd_z$ is readily seen to be extremal in $W$.
Together with $y = \pd_z$ and $h = 2z\pd_z$ it forms an $\sli_2$-triple in $W$.
Moreover, setting $v = 2z^4\pd_z$, we find
$[v,y] = 2z^3\pd_z$ and
$[v,[v,y]] = z^6\pd_z$, so $W$ is generated by $x,y,v$.
But $[y,[y,v]] = x$, so $y$ is not extremal in $W$.
\end{example}

The following result characterizes the simple Lie algebra of this example.

\begin{proposition}\label{prop_noextWitt}
Suppose that $L$ is a simple Lie algebra over the field
$\F$ of characteristic $p=5$
with an $\sli_2$-triple $x$, $y$, $h$ such that $x$ is extremal,
$-\ad_h$ is diagonizable with eigenspaces $L_i$ $(i=-2,-1,0,1,2)$ and
$[y,[y,L_{-1}]] \ne \{0\}$.
Then $L$ is isomorphic to the Witt algebra $\W$.
\end{proposition}

\begin{proof}
As $[y,[y,L_{-1}]] \subseteq L_{3} = L_{-2}$, we have
$[y,[L_{-1},y]]  = \F x$.
Let $v\in L_{-1}$ be such that $[y,[v,y]] = x$.
Consider the linear span $W$ in $L$ of $x$, $y$, $h$, $v$, $[v,y]$,
and $[v,[v,y]]$. The multiplication on these elements is fully determined:

\begin{eqnarray*}
\ [x,y] &=& h\\
\ [x,h]&=&-2x\\
\ [x,v]&=&0\ \ (\mbox{for } [x,[x,v]] \in \F x\cap L_{0} = \{0\})\\
\ [x,[v,y]] &=& [v,[x,y]]+[y,[v,x]]=[v,h]=-v\\
\ [x,[v,[v,y]]] &=& [v,[x,[v,y]]] = -[v,v]=0\\
\ [y,h] &=& 2y\\
\ [y,v] &=& -[v,y]\\
\ [y,[v,y]] &=& -[y,[y,v]] = -x\ \ \ \mbox{(by definition)}\\
\ [y,[v,[v,y]]] &=& [v,[y,[v,y]]] + 0 = [v,x]=0\\
\ [h,v]&=&v\ \ \ \mbox{(implied by the grading)}\\
\ [h,[v,y]] &=& -[v,y]\ \ \ \mbox{(implied by the grading)}\\
\ [h,[v,[v,y]]] &=& 0\ \hfill \ \mbox{(implied by the grading)}\\
\ [v,[v,y]] &=& [v,[v,y]]\\
\ [v,[v,[v,y]]] &=& 0\\
\ [[v,y],[v,[v,y]]] &=& 0
\end{eqnarray*}

Observe that $[v,[v,y]]$ is
central and that the quotient with respect to the ideal it generates
is simple of dimension $5$.  We claim
that if $[v,[v,y]] = 0$ then $W$ is isomorphic to the Witt algebra $\W$,
and otherwise $W$ is isomorphic to $\wW$, as defined in Example
\ref{ex-Witt}.

By comparison of the above multiplication rules for the spanning set
$x$, $y$, $h$, $v$, $[v,y]$, $[v,[v,y]]$ of $W$ and the basis
$-z^2 \pd_z$, $\pd_z$, $2 z \pd_z$, $2z^4 \pd_z$, $2z^3 \pd_z$, and $z^6 \pd_z$
of $\wW$ there is a surjective homomorphism
$\phi: \wW \rightarrow W$ of Lie algebras.
By assumption $W \neq 0$.
As $\wW$ only has one nontrivial proper ideal,
which maps onto $\langle [v,[v,y]] \rangle$ under $\phi$, 
it follows that $\phi$ is an isomorphism $\wW\to W$
if $[v,[v,y]] \neq 0$ and induces an isomorphism $\W \to W$ otherwise.

\begin{lettereqns}[5]{R}
It remains to prove that $L$ coincides with $W$, for then $L \cong \W$ as
$\wW$ is not simple.  To this end, suppose that $L$ strictly contains $W$, and
consider $L$ as a module on which $W$ acts.  As in the proof of Lemma
\ref{lem_ysq0}, we compute in the subalgebra $\End(L/W)$ generated by $\ad_W$.
Applying Lemma \ref{lem_ysq0}, we find that $\ad_x$ and $\ad_y$ act
quadratically on $L/\langle x,y,[x,y]\rangle$ and hence on $L/W$, so we have
relations (R1), $\ldots\,$, (R5) in $\End(L/W)$.  Write $X$, $Y$, $V$ for the
action of $\ad_x$, $\ad_y$, $\ad_v$, respectively, on $L/W$.  Due to Lemma
\ref{lem_ysq0}, and the multiplication rules $[y,[v,y]]=-x$ and $[x,[v,y]] =
-v$ listed above, we have the following relations.
\begin{eqnarray}
Y^2 &=& 0   \\
Y^2V - 2YVY + VY^2 - X &=&0  \\
XVY -XYV -VYX +YVX + V &=&0
\end{eqnarray}
Applying (R6) to (R7) and to (R2), respectively,
gives the following two relations.
\begin{eqnarray}
X + 2YVY = 0 \\
Y - YXY = 0
\end{eqnarray}
\end{lettereqns} 

Now with $R_9$, $R_{10}$ denoting the left hand sides of (R9), (R10),
respectively,
\begin{eqnarray*}
0 & = & R_9(1-XY) - 2YVR_{10} \\
  & = & (X+2YVY)(1-XY) - 2YV(Y-YXY) \\
  & = & X + 2YVY - X^2Y - 2YVYXY -2YVY + 2YVYXY \\
  & = & X.
\end{eqnarray*}
(R10) immediately implies $Y = 0$, and then (R8) implies $V=0$.

\medskip

So the images of $\ad_w$, for $w \in W$, in $\End(L/W)$ are trivial.  
This means that $W$ is an ideal of
$L$. Since $L$ is simple and $W$ is nontrivial, we find $L \cong W$,
as required.
\end{proof}

\section{The general case}\label{sec:fin}
Having dealt with the exceptional case in the previous section,
we can now proceed with the general case of Proposition
\ref{prop_Zgrading}. 

\begin{proposition}\label{prop-general}
Assume that $L$ is a simple Lie algebra over
the field $\F$ of characteristic $p\ne 2,3$,
having an $\sli_2$-pair $x$, $y$ of extremal elements.
If $L$ is not isomorphic to $\W$,
then $L$ is generated by extremal elements.
\end{proposition}

\begin{proof}
Note that $[y,[y,L_{-1}]] = 0$ as $y$ is extremal and so Proposition
\ref{prop_Zgrading} gives that $h = [x,y]$ is diagonizable and the
components $L_i = L_i(-\ad_h)$ $(i=-2,-1,0,1,2)$ of the grading by $h$
satisfy $L_{-2} = \F x$, $L_{-1} = [x,L_1]$, $L_2=\F y$,
and $L_1=[y,L_{-1}]$.

Consider the subalgebra $I$ of $L$ generated by
$x$, $y$, and $L_{1}$. As $L_{-1}=[x,L_1]$,
the subalgebra $I$ contains the linear
subspace $J=L_{-2}+L_{-1}+L_1+L_2$ of $L$.
As $[J,L_0]\subseteq J$ and $J$ generates $I$, 
we have $[I,L_0]\subseteq I$. This
implies
$[I,L]=I$. In other words, $I$ is an ideal of $L$, and so,
by simplicity of $L$, it coincides with $L$. 
Therefore, it suffices to show that
for each $z\in L_1$ there exists an
extremal element $u\in L$
such that $z$ is in the subalgebra
generated by $x$, $y$, and $u$.

To this end, let $z\in L_1$.
Put $h=[x,y]$. The following
relations hold in $L$, for some $\alpha\in\F$.
\begin{eqnarray}
\label{rel1}
[h,x]=2x, && \\
\label{rel2}
[h,y]=-2y, && \\
\label{rel3}
[z,h]=z, && \\
\label{rel4}
[y,z]=0, && \\
\label{rel5}
[x,[x,z]]=0, && \\
\label{rel6}
[y,[x,z]]=z, && \\
\label{rel7}
[y,[z,[z,x]]]=0, && \\
\label{rel8}
[x,[z,[z,x]]]=0, && \\
\label{rel9}
[y,[z,[z,[z,x]]]]=0,&& \\ 
\label{rel10}
[x,[x,[z,[z,[z,x]]]]]=0, && \\
\label{rel11}
[y,[x,[z,[z,[z,x]]]]]=[z,[z,[z,x]]], && \\
\label{rel12}
[z,[z,[z,[z,x]]]] = \alpha  y. &&
\end{eqnarray}

We claim that the Lie subalgebra $L'$ of $L$ generated by $x$, $y$, and $z$
is linearly spanned by the following set $B$ of eight elements, where $h_1 =
[[x,z],z]$.
\[
x\in L_{-2};\ 
\, [x,z], [[h_1,z],x]\in L_{-1};\ 
h, [[x,z],z]]\in L_0 ;\ 
z, [h_1,z]\in L_1 ;\ 
y\in L_2
\]
To see that this is true, we verify that the images of the elements of $B$
under $\ad_x$, $\ad_y$, and $\ad_z$ are scalar multiples of these.  For
$\ad_x$ and $\ad_y$ this is straightforward.  
As for $\ad_z$, the statement is trivially verified for
all elements of $B$ but
$[[h_1,z],x]$.
As $[h_1,x]=0$, we have
\begin{eqnarray*}
\ad_z([[h_1,z],x])
&=&[\ad_z([h_1,z]),x] + 
[[h_1,z],\ad_z(x)]\\
&=&\alpha [y,x] + [[h_1,\ad_z(x)],z]
+[h_1,[z,\ad_z(x)]]\\
&=&-\alpha h -\ad_z([h_1,\ad_z(x)]) + [h_1,h_1] \\
&=&-\alpha h -\ad_z([[h_1,z],x]), 
\end{eqnarray*}
so $\ad_z([[h_1,z],x])= -\frac{\alpha}{2}h$.
This establishes the claim that $L'$ is linearly spanned by $B$.

We exhibit an element $u\in L'$ as specified.
Because of the 
grading induced by $\ad_h$ on $L$, the endomorphism
$\ad_z$ on $L$ is nilpotent of
order at most $5$ and
$\exp(\ad_z)$ 
is a linear transformation of $L$
(it is well defined as $p\ne2,3$).
Put $$u = \exp(\ad_z)x
= x +\ad_z(x)+\frac{1}{2}\ad_z^2(x)+\frac{1}{6}\ad_z^3(x)
+\frac{1}{24}\ad_z^4(x).$$ 
A straightforward computation in $L'$ shows
that $y$ and $u$ are an $\sli_2$-pair in $L$.
By (\ref{rel7}), (\ref{rel9}), (\ref{rel12}), and (\ref{rel6}) we find

\begin{eqnarray}\label{relyu}
[y,u] & = & [y, x] + [y,[z,x]] + \frac{1}{2}[y,\ad_z^2(x)] + \frac{1}{6}[y,\ad_z^3(x)] + \frac{1}{24}[y,\ad_z^4(x)], \cr
& = & [y,x] + [y,[z,x]] + 0 + 0 + [y,\alpha y] \cr
& = & -h -z
\end{eqnarray}
so, by (\ref{rel2}), (\ref{rel4}),
\begin{equation}\label{eqn_yuy}
[[y,u],y] =  -[h,y] -[z,y] = 2y.
\end{equation}

For $[[y,u],u]$ we compute, using (\ref{relyu}), (\ref{rel1}), and Proposition \ref{prop_Zgrading},
\begin{eqnarray*}
[[y,u],x] & = & -[h,x] -[z,x] = -2x -[z,x], \cr
[[y,u], \ad_z(x)] & = & -\ad_z(x) - \ad_z^2(x), \cr
[[y,u], \ad_z^2(x)] & = & 0 -\ad_z^3(x), \cr
[[y,u], \ad_z^3(x)] & = & \ad_z^3(x) - \ad_z^4(x), \cr
[[y,u], \ad_z^4(x)] & = & +2\ad_z^4(x) -\ad_z^5(x) = 2\ad_z^4(x),
\end{eqnarray*}
so
\begin{eqnarray}\label{eqn_yuu}
[[y,u],u] & = & -2x -[z,x]-\ad_z(x) - \ad_z^2(x) - \frac{1}{2}\ad_z^3(x) + \cr
& & \frac{1}{6}(\ad_z^3(x) - \ad_z^4(x)) + \frac{1}{24}(2 \ad_z^4(x)) \cr
& = & -2u,
\end{eqnarray}

Now (\ref{eqn_yuy}) and (\ref{eqn_yuu}) show that $y$ and $u$ are an $\sli_2$-pair in $L$.

\bigskip

By Propositions \ref{prop_Zgrading} and \ref{prop_noextWitt},
and the assumption that $L$ is not isomorphic to $\W$,
this implies that $u$ is extremal in $L$.

We verify that $z$ lies in the subalgebra $L''$ of $L$
generated by the three
extremal elements $x$, $y$, and $u$.
Observe that
$$\ad_z(x)+\frac{1}{2}\ad_z^2(x)+\frac{1}{6}\ad_z^3(x)
= u - x -\frac{\alpha}{24}y\in L''.$$ 
Acting by $\ad_y$ and using (\ref{rel6}), (\ref{rel7}), (\ref{rel9}), we find
\[
 z = - \ad_y\ad_z(x)-\frac{1}{2}\ad_y\ad_z^2(x)-\frac{1}{6}\ad_y \ad_z^3(x)
\in\ad_y L''\subseteq L''.
\]
This proves that $z$ belongs to $L''$ and so we are done.
\end{proof}

\begin{proof}[Proof of Theorem \ref{th:main}]
Let $L$ be as in the assumption.
By Theorem \ref{thm_sl2quad}, there is an $\sli_2$-pair $x,y$ in $L$
with $x$ extremal in $L$. 
Proposition \ref{prop-general} finishes the proof.
\end{proof}

\bibliographystyle{alpha}

\end{document}